\documentclass[11pt]{amsart}
\usepackage[mathscr]{eucal}
\usepackage{amssymb, amsmath,array, amscd}
\usepackage{enumerate,verbatim}
\usepackage{graphicx}
\usepackage{url}
\usepackage[plainpages]{hyperref}

\input xy
\xyoption{all}

\newtheorem{theorem}{Theorem}

\theoremstyle{definition}

\newtheorem{remark}[theorem]{Remark}

\newcommand{\set}[1]{\left\{#1\right\}}

\renewcommand{\k}{\Bbbk}

\newcommand{\RR}{{R}}
\newcommand{\bC}{{C}}

\newcommand{\R}{\mathbb{R}}
\newcommand{\Z}{\mathbb{Z}}

\newcommand{\CC}{{\mathcal C}}

\newcommand{\PS}{{P\!\varSigma}}

\DeclareMathOperator{\Hom}{Hom}

\DeclareMathOperator{\Span}{span}

\begin{document}

\title{Resonance of basis-conjugating automorphism groups}

\author[Daniel C. Cohen]{Daniel C. Cohen$^\dag$}
\address{Department of Mathematics, Louisiana State University,
Baton Rouge, LA 70803, USA}
\email{\href{mailto:cohen@math.lsu.edu}{cohen@math.lsu.edu}}
\urladdr{\href{http://www.math.lsu.edu/~cohen/}
{http://www.math.lsu.edu/\char'176cohen}}
\thanks{{$^\dag$}Partially supported 
by National Security Agency grant H98230-05-1-0055}

%14M12 Determinantal varieties
%20F28 Automorphism groups of groups
%20J06 Cohomology of groups


\begin{abstract}
We determine the structure of the first resonance variety of the cohomology ring of the group of automorphisms of a finitely generated free group which act by conjugation on a given basis.
\end{abstract}

%\date{\today}

\maketitle

\section{Resonance of $\PS_n$} \label{sec:intro}

Let $F_n$ be the free group generated by $x_1,\dots,x_n$.  The basis-conjugating automorphism group, or pure symmetric automorphism group, is the group $\PS_n$ of all automorphisms of $F_n$ which send each generator $x_i$ to a conjugate of itself.  Results of Dahm \cite{Dahm} and 
Goldsmith \cite{Goldsmith} imply that this group may also be realized as the ``group of loops,'' 
the group of motions of a collection of $n$ unknotted, unlinked oriented circles in $3$-space, where each circle returns to its original position.  McCool \cite{McC} found the following presentation for the basis-conjugating automorphism group:
\begin{equation}\label{eq:McCool}
\PS_n = \langle \beta_{i,j}, 1\le i \neq j \le n \mid [\beta_{i,j},\beta_{k,l}], [\beta_{i,k},\beta_{j,k}], [\beta_{i,j},(\beta_{i,k}\cdot \beta_{j,k})]
\rangle,
\end{equation}
where $[u,v]=u^{}v^{}u^{-1}v^{-1}$ denotes the commutator, the indices in the relations are distinct, and the generators $\beta_{i,j}$ are the automorphisms of $F_n$ defined by
\[
\beta_{i,j}(x_k)=\begin{cases} x_k &\text{if $k \neq j$,}\\ x_j^{-1} x_i^{} x_j^{}&\text{if $k=i$.}
\end{cases}
\]
The purpose of this note is to determine the structure of the first resonance variety of the cohomology ring of this group.

Let $A=\bigoplus_{k=0}^\ell A^k$ be a finite-dimensional, graded, connected 
algebra over an algebraically closed field $\k$ of 
characteristic $0$.   Since $a\cdot a=0$ for each $a\in A^1$, multiplication 
by $a$ defines a cochain complex $(A,\delta_a)$:
\[
A^0\xrightarrow{\quad \delta_a\quad} A^1 \xrightarrow{\quad \delta_a\quad} A^2 \xrightarrow{\quad \delta_a\quad} \cdots\cdots \xrightarrow{\quad \delta_a\quad} A^\ell,
\]
where $\delta_a(x)=ax$.  
The resonance varieties of $A$ are the 
jumping loci for the cohomology of these complexes:
$R^{j}_{d}(A)=\set{a\in A^1 \mid \dim_{\k} H^j(A,\delta_a) \ge d}$.  As shown by Falk \cite{Falk}, these algebraic subvarieties of $A^1$ are isomorphism-type invariants of the algebra $A$.  
They have 
been the subject of considerable recent interest, particularly 
in the context of hyperplane arrangements, see for instance Yuzvinsky 
\cite{Yuz}, and references therein.
We will focus on the first resonance variety $R^1(A)=\set{a\in A^1 \mid H^1(A,\delta_a)\neq 0}$.  

Since the relations in 
the presentation \eqref{eq:McCool} of $\PS_n$ 
are all commutators, the first homology group $H_1(\PS_n;\k)$ is a vector space of dimension $n(n-1)$ with basis $\set{[\beta_{p,q}] \mid 1\le p \neq q \le \ell}$.  Let $\set{e_{p,q} \mid 1\le p \neq q \le \ell}$ be the dual basis of $H^1(\PS_n;\k)$.  
Denote the first resonance variety of $A=H^*(\PS_n;\k)$ by $R^1(\PS_n,\k)$.

\begin{theorem} \label{thm:res}
The first resonance variety of the cohomology ring $A=H^*(\PS_n;\k)$ of the basis-conjugating autormorphism group is 
\[
R^1(\PS_n,\k) = \bigcup_{1\le i < j \le n}\!\! \Span\set{e_{i,j},e_{j,i}} \ \cup 
\bigcup_{1\le i < j<k \le n}\!\!\! \Span\set{e_{j,i}-e_{k,i},e_{i,j}-e_{k,j},e_{i,k}-e_{j,k}}.
\]
\end{theorem}

This result reveals an interesting relationship between the resonance variety and the Bieri-Neumann-Strebel (BNS) invariant of the basis-conjugating automorphism group.  For a finitely generated group $G$, let $\CC$ be the Cayley graph corresponding to a finite generating set.  Given an additive character $\chi\colon G \to \R$, let $\CC_+(\chi)$ be the full subgraph of $\CC$ on vertex set $\set{g\in G \mid \chi(g) \ge 0}$.  Then the (first) BNS invariant of $G$ is the conical subset $\Sigma(G)$ of $\Hom(G,\R)=H^1(G;\R)$ defined by
\[
\Sigma(G)=\set{\chi \in \Hom(G,\R) \setminus\set{0} \mid \CC_+(\chi)\ \text{is connected}}.
\]
This invariant of $G$ (which is independent of the choice of generating set) may be used 
to determine which subgroups above the commutator subgroup $[G,G]$ are finitely generated, see \cite{BNS}.

The BNS invariant of the group $G=\PS_n$ was determined by Orlandi-Korner \cite{O-K}.  Combining her result with the above theorem yields the following.

\begin{theorem} \label{cor:BNS}
The Bieri-Neumann-Strebel invariant of the basis-conjugating automorphism group is given by
\[
\Sigma(\PS_n) = H^1(\PS_n;\R) \setminus R^1(\PS_n,\R).
\]
\end{theorem}
This relationship between the resonance variety and BNS invariant is known to hold for some other groups, including right-angled Artin groups, see Meier and VanWyk~\cite{MVW} and Papadima and Suciu \cite{PS06}.

\section{Proof of Theorem \ref{thm:res}}
\label{sec:PSn}

In this section, we recall the structure of the cohomology ring of the basis-conjugating automorphism group, and use it to prove Theorem \ref{thm:res}.  The cohomology of $\PS_n$ was 
computed by Jensen, McCammond, and Meier \cite{JMM}, resolving positively a conjecture of Brownstein and Lee \cite{BL}.

\begin{theorem}[{\cite{JMM}}] \label{thm:cohomology}
Let $E_\Z$ denote the exterior algebra over $\Z$ generated by degree one elements $e_{p,q}$, $1\le p \neq q \le n$, and let $I_\Z$ denote the two-sided ideal in $E_\Z$ generated by 
\begin{equation} \label{eq:Igens}
\begin{aligned} 
\eta_{i,j} &= e_{i,j} e_{j,i},\  1\le i< j\le n, \hfill \\
\tau_{i,j}^k &= (e_{k,i}-e_{j,i})(e_{k,j}-e_{i,j}),
\  1\le k \le n,\  1\le i<j\le n,\  k\notin\set{i,j}.
\end{aligned}
\end{equation}
Then the integral cohomology algebra of the basis-conjugating automorphism group $\PS_n$ is isomorphic to the quotient of $E$ by $I$, $H^*(\PS_n;\Z) \cong E_\Z/I_\Z$.  
\end{theorem}

\begin{remark}
The above presentation of the cohomology ring $H^*(\PS_n;\Z)$ differs slightly from that given in \cite{BL,JMM}, but is easily seen to be equivalent.  For instance, the relation in $H^*(\PS_n;\Z)$ arising from the generator $\tau^k_{i,j}$ of the ideal $I$ may be obtained by constructing an appropriate linear combination of the relations labeled $2$ and $3$ in \cite[Thm. 6.7]{JMM}. 
\end{remark}

Let $\k$ be a field of characteristic zero.  From Theorem \ref{thm:cohomology} and the Universal Coefficient Theorem, the cohomology algebra $H^*(\PS_n;\k)$ is isomorphic to $E_\k/I_\k$, where $E_\k$ is the exterior algebra over $\k$  generated by $e_{p,q}$, $1\le p\neq q\le n$, 
and $I_\k$ is the ideal in $E_\k$ generated by the elements $\eta_{i,j}$ and $\tau^k_{i,j}$ above.  

Recall from the Section \ref{sec:intro} that the first resonance variety of $A=H^*(\PS_n;\k)$ is
\[
R^1(\PS_n,\k) = \{ a \in A^1 \mid H^1(A,\delta_a) \neq 0\}.
\]
Observe that $A^1=H^1(\PS_n;\k)$ is a vector space of dimension $N=n(n-1)$ over $\k$.  Elements of $A^1$ are of the form $a = \sum_{p\neq q} a_{p,q} e_{p,q}$, where $a_{p,q} \in \k$.  For $1\le i < j \le n$, write
\begin{equation} \label{eq:Cij}
C_{i,j} = \Span\set{e_{i,j},e_{j,i}}=\{ a \in A^1 \mid a_{p,q}=0\text{ if } \{p,q\} \neq \{i,j\}\},
\end{equation}
and for $1\le i<j<k \le n$, write
\begin{equation} \label{eq:Cijk}
\begin{aligned} 
C_{i,j,k}&=\Span\set{e_{j,i}-e_{k,i},e_{i,j}-e_{k,j},e_{i,k}-e_{j,k}}\\
&=
\Biggl\{a \in A^1 \Biggm|
\begin{matrix}
a_{j,i}+a_{k,i}=0,\ a_{i,j}+a_{k,j}=0,\ a_{i,k}+a_{j,k}=0,\\[2pt]
a_{p,q}=0\text{ if } \{p,q\} \not\subset \{i,j,k\} \hfill \  
\end{matrix}
\Biggr\}.
\end{aligned}
\end{equation}
Note that these are linear subspaces of $A^1\cong \k^N$ of dimensions $2$ and $3$ respectively. 
In this notation, Theorem \ref{thm:res} asserts that $R^1(\PS_n,\k)=\bigcup_{1\le i<j\le n} C_{i,j} \cup \bigcup_{1\le i<j<k\le n} C_{i,j,k}$.

\begin{proof}[Proof of Theorem \ref{thm:res}]
In the case $n=2$, we have $\PS_2=F_2=\langle \beta_{1,2},\beta_{2,1}\rangle$, and the theorem asserts that $R^1(\PS_2,\k) = C_{1,2} = H^1(\PS_2;\k)$, which is clear.  So assume that $n \ge 3$.

Write $\RR=R^1(\PS_n,\k)$ and $\bC=\bigcup_{1\le i < j \le n} C_{i,j} \cup 
\bigcup_{1\le i < j<k \le n} C_{i,j,k}$.  Observe that $0 \in \bC$ and $0 \in \RR$.
So it is enough to show that $\RR \setminus\{0\}=\bC\setminus\{0\}$.

Write $E=E_\k$ and $I=I_\k$.  
For $a \in A^1 = E^1$, we have a short exact sequence of chain complexes
$0 \longrightarrow (I,\delta_a) \xrightarrow{\ \iota\ } (E,\delta_a) \xrightarrow{\ p\ } (A,\delta_a) \longrightarrow 0$:
\begin{equation*}
\xymatrixrowsep{20pt}
\xymatrixcolsep{40pt}
\xymatrix{
&& I^2 \ar[d]^{\iota^2} \ar[r]^{\delta_a} 
& I^3 \ar[d]^{\iota^3} \ar[r]
& \dots 
\\
E^0 \ar[d]^{p^0}  \ar[r]^{\delta_a}
& E^1 \ar[d]^{p^1} \ar[r]^{\delta_a}
& E^2 \ar[d]^{p^2} \ar[r]^{\delta_a}
& E^3 \ar[d]^{p^3} \ar[r] & \dots
 \\
A^0 \ar[r]^{\delta_a}
& A^1 \ar[r]^{\delta_a}
& A^2 \ar[r]^{\delta_a}
& A^3  \ar[r] & \dots
}
\end{equation*}
where $\iota\colon I \to E$ is the inclusion, $p\colon E \to A$ the projection, and $\delta_a(x) = a x$.  Note that, since $I$ is generated in degree two,  the maps $p^0\colon E^0 \to A^0$ and $p^1\colon E^1 \to A^1$ are identity maps.  If $a \neq 0$, then the complex $(E,\delta_a)$ is acyclic.  Consequently, the corresponding long exact cohomology sequence yields
\[
H^1(A,\delta_a) \cong H^2(I,\delta_a) = \ker(\delta_a\colon I^2 \to I^3) = 
\ker(\iota^3 \circ \delta_a\colon I^2 \to E^3).
\]
Thus, $a\in \RR\setminus\{0\}$ if and only if the map $\psi_a:=\iota^3 \circ \delta_a$ fails to inject.

Since the elements $\eta_{i,j}$ and $\tau^k_{i,j}$ recorded in \eqref{eq:Igens} generate the ideal $I$ and are of degree two (and are linearly independent in $E^2$), these elements form a basis for $I^2$.  We record the images of these basis elements under the map $\psi_a$.  For $1\le i<j\le n$, 
\begin{equation} \label{eq:eta}
\psi_a(\eta_{i,j}) =  \sum_{\{p,q\}\neq\{i,j\}} a_{p,q} e_{p,q} \eta_{i,j}
= \sum_{\{p,q\}\neq\{i,j\}} a_{p,q} e_{p,q} e_{i,j} e_{j,i}.
\end{equation}
For $1\le k \le n$, $1\le i <j\le n$, $k \notin\set{i,j}$,  
\begin{equation} \label{eq:tau}
\begin{aligned} 
\psi_a(\tau^k_{i,j})&= (a_{j,i}+a_{k,i})e_{j,i}e_{k,i}(e_{k,j}-e_{i,j})
-(a_{i,j}+a_{k,j})e_{i,j}e_{k,j}(e_{k,i}-e_{j,i})  \\
&\qquad+(a_{i,k}e_{i,k}+a_{j,k}e_{j,k}) \tau^k_{i,j}
+\sum_{\{p,q\}\not\subset\{i,j,k\}} a_{p,q}e_{p,q} \tau^k_{i,j}.
\end{aligned}
\end{equation}

These calculations immediately yield the containment $\bC\setminus\{0\} \subseteq \RR \setminus\{0\}$.  If $a \in C_{i,j}$, then $a_{p,q}=0$ for $\{p,q\}\neq \{i,j\}$.  For such $a$, we have $\psi_a(\eta_{i,j})=0$ by \eqref{eq:eta}, so $C_{i,j} \subset \RR$.
If $1\le i<j<k\le n$ and $a \in C_{i,j,k}$, then  $a_{j,i}+a_{k,i}=0$, $a_{i,j}+a_{k,j}=0$, $a_{i,k}+a_{j,k}=0$, and 
$a_{p,q}=0$ for $\{p,q\} \not\subset \{i,j,k\}$.  In this instance, \eqref{eq:tau} may be used to check that $a_{j,i} \tau^k_{i,j}-a_{i,k}\tau^i_{j,k}$ and 
$a_{i,j} \tau^k_{i,j}-a_{i,k}\tau^j_{i,k}$ are elements of $\ker(\psi_a)$.  If $a \in C_{i,j,k}$ is nonzero, at least one of $a_{i,j},a_{i,k},a_{j,i}$ must be nonzero.
Consequently, $\psi_a$ has nontrivial kernel, and $C_{i,j,k} \subset \RR$.

Establishing the reverse containment, $\RR \setminus \{0\} \subseteq \bC \setminus \{0\}$, is more involved.  We will show that $a \notin \bC$ implies 
that $a \notin \RR$.  If $a \notin \bC$, then $a \neq 0$.  So assume without loss that $a_{2,1} \neq 0$.  Since $a \notin C_{1,2} \subset \bC$, we must also have 
$a_{p,q} \neq 0$ for some $\{p,q\} \neq \{1,2\}$.  We will consider several cases depending on the relationship between the sets $\{1,2\}$ and $\{p,q\}$.

\smallskip

\noindent\textbf{Case 1.} \quad $\set{1,2} \cap \set{p,q}=\emptyset$\\
Assume first that $\{1,2\}$ and $\{p,q\}$ are disjoint. Note that $n\ge 4$ in this instance. 
Permuting indices if necessary, we may assume that $a\in H^1(\PS_n;\k)$ satisfies $a_{2,1}\neq 0$ and $a_{3,4}\neq 0$. We will show that this assumption implies that the map $\psi_a\colon I^2 \to E^3$ injects, hence $a \notin R$.  Specifically, we will exhibit a subspace $V \subset 
E^3$ and a projection $\pi\colon E^3\twoheadrightarrow V$ so that the composition $\pi \circ \psi_a\colon I^2 \to V$ is an isomorphism. 

Let $\mathcal{V}$ be the union of the sets
\[
\begin{aligned}
&\set{e_{1,2}e_{2,1}e_{3,4}}\cup\set{e_{2,1} e_{i,j} e_{j,i} \mid 1\le i<j \le n,\ \set{i,j}\neq \set{1,2}},\\
&\set{e_{3,4}e_{1,2}e_{k,1},\ e_{3,4}e_{2,1}e_{1,k},\ e_{3,4}e_{1,2}e_{1,k}\mid 3\le k\le n},\\
&\set{e_{2,1}e_{k,i}e_{k,j},\ e_{2,1}e_{j,i}e_{j,k},\ e_{2,1}e_{i,k}e_{k,j}\mid 1\le i \le 2<j< k\le n
\ \text{or}\ 3\le i<j<k\le n}. 
\end{aligned}
\]
There is a bijection between $\mathcal{V}$ and the set of generators of $I^2$ given by
\[
\begin{matrix}
\eta_{1,2}\leftrightarrow e_{1,2}e_{2,1}e_{3,4},&\eta_{i,j} \leftrightarrow e_{2,1}e_{i,j}e_{j,i}&(\set{i,j}\neq\set{1,2}),\\
\tau^k_{1,2} \leftrightarrow e_{3,4}e_{1,2}e_{k,1},& \tau^2_{1,k} \leftrightarrow e_{3,4}e_{2,1}e_{1,k},& \tau^1_{2,k} \leftrightarrow e_{3,4}e_{1,2}e_{1,k} &(3\le k\le n),\hfill \ \\
\tau^k_{i,j} \leftrightarrow e_{2,1}e_{k,i}e_{k,j},& \tau^j_{i,k} \leftrightarrow e_{2,1}e_{j,i}e_{j,k},& \tau^i_{j,k} \leftrightarrow e_{2,1}e_{i,k}e_{k,j} &(1\le i\le 2<j< k\le n).
\end{matrix}
\]
In particular, the monomials in $\mathcal{V}$ are distinct. Hence, $\mathcal{V}$ is a linearly independent set in $E^3$ of cardinality $|\mathcal{V}|=\dim_\k I^2=\binom{n}{2}(n-1)$.
Let $V=\Span\mathcal{V}\subset E^3$.

Define $\pi\colon E^3 \twoheadrightarrow V$ on basis elements by $\pi(e_{a,b}e_{c,d}e_{p,q}) = e_{a,b}e_{c,d}e_{p,q}$ if 
$e_{a,b}e_{c,d}e_{p,q} \in \mathcal{V}$, and $\pi(e_{a,b}e_{c,d}e_{p,q}) = 0$ if 
$e_{a,b}e_{c,d}e_{p,q} \notin \mathcal{V}$.  Then, a calculation using 
\eqref{eq:eta} and  \eqref{eq:tau} reveals that $\pi \circ \psi_a\colon I^2 \to V$ is an isomorphism.  For instance, ordering the bases of $I^2$ and $V$ appropriately,  one can check that the matrix $\sf{M}$ of 
$\pi \circ \psi_a$ has determinant $\det {\sf{M}}=y_{2,1}^{m_{2,1}} y_{3,4}^{m_{3,4}}\neq 0$, where $m_{3,4}=3n-5$ and $m_{2,1}=\binom{n}{2}(n-1)-m_{3,4}$. 
Thus, if $a\notin C$ satisfies $a_{2,1}\neq 0$ and $a_{p,q}\neq 0$ for some $p,q$ with $\set{1,2}\cap\set{p,q}=\emptyset$, then $a \notin R$.

\smallskip

\noindent\textbf{Case 2.} \quad $\set{1,2} \cap \set{p,q}\neq\emptyset$\\
Now assume that $a\notin C$, $a_{2,1}\neq 0$, and $a_{r,s}=0$ for all $r,s$ with $\set{1,2}\cap\set{r,s}=\emptyset$.  Since $a\notin C_{1,2}\subset C$, we must have $a_{p,q}\neq 0$ for some $p,q$ with $|\{p,q\}\cap\{1,2\}|=1$.  Permuting indices if necessary, we may assume that $3 \in \{p,q\}$.

In the case $n=3$, since $a \notin C_{1,2,3}\subset C$, one of the sums $a_{2,1}+a_{3,1}$, $a_{1,2}+a_{3,2}$, $a_{1,3}+a_{2,3}$ must be nonzero, see \eqref{eq:Cijk}.  
In this instance, ordering bases appropriately, the map $\psi_a\colon I^2 \to E^3$ has matrix
\begin{equation} \label{eq:M3}
{\sf{M}}_3=
{\footnotesize{
\left(
\begin{array}{cccccc}
 a_{3,2} & 0 & 0 & 0 & -a_{1,2}-a_{3,2} & 0 \\
 a_{3,1} & 0 & 0 & 0 & -a_{2,1}-a_{3,1} & 0 \\
 a_{2,3} & 0 & 0 & a_{1,2} & -a_{2,3} & a_{2,1} \\
 -a_{1,3} & 0 & 0 & a_{1,2} & a_{1,3} & a_{2,1} \\
 0 & a_{3,2} & 0 & -a_{3,2} & a_{1,3} & -a_{3,1} \\
 0 & -a_{2,3} & 0 & a_{1,3}+a_{2,3} & 0 & 0 \\
 0 & -a_{2,1} & 0 & a_{2,1}+a_{3,1} & 0 & 0 \\
 0 & a_{1,2} & 0 & -a_{1,2} & -a_{1,3} & a_{3,1} \\
 0 & 0 & -a_{3,1} & a_{3,2} & a_{2,3} & a_{3,1} \\
 0 & 0 & a_{2,1} & a_{3,2} & a_{2,3} & -a_{2,1} \\
 0 & 0 & a_{1,3} & 0 & 0 & -a_{1,3}-a_{2,3} \\
 0 & 0 & a_{1,2} & 0 & 0 & -a_{1,2}-a_{3,2} \\
 0 & 0 & 0 & a_{3,2} & -a_{1,3} & -a_{2,1} \\
 0 & 0 & 0 & a_{2,1}+a_{3,1} & 0 & 0 \\
 0 & 0 & 0 & a_{1,3}+a_{2,3} & 0 & 0 \\
 0 & 0 & 0 & a_{1,2} & -a_{2,3} & -a_{3,1} \\
 0 & 0 & 0 & 0 & a_{1,2}+a_{3,2} & 0 \\
 0 & 0 & 0 & 0 & a_{2,1}+a_{3,1} & 0 \\
 0 & 0 & 0 & 0 & 0 & a_{1,2}+a_{3,2} \\
 0 & 0 & 0 & 0 & 0 & a_{1,3}+a_{2,3}
\end{array}
\right)}}.
\end{equation}
Using the assumptions on $a_{2,1}$,  $a_{2,1}+a_{3,1}$, $a_{1,2}+a_{3,2}$, $a_{1,3}+a_{2,3}$, and $a_{p,q}\in\set{a_{1,3},a_{2,3},a_{3,1},a_{3,2}}$, it is readily checked that the matrix ${\sf{M}}_3$ has maximal rank.  Hence, if $a\notin C$, then 
$\psi_a\colon I^2 \to E^3$ injects in the case $n=3$.

For general $n$, the assumption that $a\notin C_{1,2,3}$ implies that the set
\[
\set{a_{2,1}+a_{3,1}, a_{1,2}+a_{3,2}, a_{1,3}+a_{2,3}}\cup\set{a_{r,s} \mid \set{r,s}\not\subset\set{1,2,3}}
\] 
contains a nonzero element.  Recall that, by Case 1, we may assume that $a_{r,s}=0$ for all $r,s$ with $\set{1,2}\cap \set{r,s}=\emptyset$. 
If $a_{r,s}=0$ for all $\set{r,s}\not\subset\set{1,2,3}$, let $W$ be the subspace of $E^3$ spanned by the union of the sets
\[
\begin{aligned}
&\set{e_{i_1,j_1}e_{i_2,j_2}e_{i_3,j_3}\mid 1\le i_k,j_k\le 3,\ i_k\neq j_k},\\
&\set{e_{1,2}e_{2,1}e_{p,q}}\cup\set{e_{2,1} e_{i,j} e_{j,i} \mid 1\le i<j \le n,\ \set{i,j}\neq \set{1,2}},\\
&\set{e_{2,1}e_{k,i}e_{k,j},\ e_{2,1}e_{j,i}e_{j,k},\ e_{2,1}e_{i,k}e_{k,j}\mid 1\le i \le 2<j< k\le n
\ \text{or}\ 3\le i<j<k\le n},\\
&\set{
\begin{matrix}
e_{p,q}e_{k,1}e_{k,2},&e_{p,q}e_{k,1}e_{1,2},&e_{p,q}e_{2,1}e_{k,2},\\
e_{p,q}e_{2,1}e_{2,k},&e_{p,q}e_{2,1}e_{1,k},&e_{p,q}e_{k,1}e_{2,k},\\
e_{p,q}e_{1,2}e_{1,k},&e_{p,q}e_{1,2}e_{2,k},&e_{p,q}e_{k,2}e_{1,k}
\end{matrix}
\Biggm| 3\le k \le n
}.
\end{aligned}
\]
Define $\pi\colon E^3 \twoheadrightarrow W$ on basis elements as before.  Ordering bases appropriately, one can use \eqref{eq:eta} and \eqref{eq:tau} to find a submatrix $\sf{M}$ of the matrix of $\pi\circ \psi_a\colon I^2 \to W$ of the form
\[
{\sf{M}}=\left(
\begin{matrix}
{\sf{U}} & * \\
0 & {\sf{M}}_3
\end{matrix}
\right),
\]
where ${\sf{M}}_3$ is given by \eqref{eq:M3} and $\sf{U}$ is upper triangular, with diagonal entries $a_{2,1}\neq 0$ and $a_{p,q}\neq 0$.  
(The choice of $\sf{U}$ depends on which $a_{p,q}\in\set{a_{1,3},a_{2,3},a_{3,1},a_{3,2}}$ is nonzero.) 
Consequently, the matrix $\sf{M}$ has maximal rank. It follows that 
$\psi_a\colon I^2 \to E^3$ injects.

Finally, consider the case where $a\notin C_{1,2,3} \subset C$, $a_{2,1}\neq 0$, 
$a_{p,q}\neq 0$ for some $a_{p,q}\in\set{a_{1,3},a_{2,3},a_{3,1},a_{3,2}}$, 
and $a_{r,s}\neq 0$ for some $\set{r,s}\not\subset\set{1,2,3}$.  
Since we may assume by Case 1 that $a_{r,s}=0$ if $\set{1,2}\cap\set{r,s}=\emptyset$, we have $a_{r,s}\neq 0$ for some $r,s$ with $r\in\set{1,2}$ and $4\le s \le n$. 
In this instance, let $W$ be the subspace of $E^3$ spanned by the union of the sets
\[
\begin{aligned}
&\set{e_{1,2}e_{2,1}e_{p,q}}\cup\set{e_{2,1} e_{i,j} e_{j,i} \mid 1\le i<j \le n,\ \set{i,j}\neq \set{1,2}},\\
&\set{e_{2,1}e_{k,i}e_{k,j},\ e_{2,1}e_{j,i}e_{j,k},\ e_{2,1}e_{i,k}e_{k,j}\mid 1\le i \le 2<j< k\le n
\ \text{or}\ 3\le i<j<k\le n},\\
&\set{e_{r,s}e_{3,1}e_{3,2},\ e_{r,s}e_{1,3}e_{3,2},\ e_{r,s}e_{2,3}e_{3,1}},\\
&\set{e_{p,q}e_{k,1}e_{k,2},\ e_{p,q}e_{2,k}e_{k,1},\ e_{p,q}e_{1,k}e_{k,2} \mid 4\le k \le n}.
\end{aligned}
\]
Defining $\pi\colon E^3 \twoheadrightarrow W$ on basis elements as above, a calculation using 
\eqref{eq:eta} and \eqref{eq:tau} shows that $\pi \circ \psi_a\colon I^2 \to W$ is an isomorphism.  For instance, ordering bases appropriately, one can check that the matrix $\sf{M}$ of $\pi \circ \psi_a$ has determinant $\det{\sf{M}}=a_{2,1}^{m_{2,1}}a_{p,q}^{m_{p,q}}a_{r,s}^{3}$, where 
$m_{p,q}=3n-8$ and $m_{2,1}=\binom{n}{2}(n-1)-m_{p,q}$.  Hence, $\psi_a$ injects in this final case.

Thus, for any $a\notin C$, the map $\psi_a\colon I^2 \to E^3$ injects, and $a \notin R$.  This completes the proof of Theorem \ref{thm:res}.
\end{proof}

\begin{remark}
It follows from Theorem \ref{thm:cohomology} that the integral cohomology groups of $\PS_n$ are torsion free, with Betti numbers $b_k(\PS_n)=\text{rank}\,H^k(\PS_n;\Z)$ given by the coefficients of the 
Poincar\'e polynomial  ${\mathfrak{p}}(\PS_n,t) = \sum_{k\ge 0} b_k(\PS_n)\cdot t^k=
(1+nt)^{n-1}$, see \cite[\S{6}]{JMM}.  Thus the cohomology groups cannot distinguish $\PS_n$ from a direct product 
$F_n\times\dots\times F_n$ of $n-1$ free groups of rank $n$.

These groups are, however, distinguished by their cohomology rings for $n\ge 3$.  By Theorem \ref{thm:res}, the irreducible components of $R^1(\PS_n,\k)$ are two- and three-dimensional.  On the other hand, the results of \cite{CScv} or \cite{PS06} may be used to show that the irreducible components of the first resonance variety of $H^*(F_n\times\dots\times F_n;\k)$ are all $n$-dimensional.  
\end{remark}

\section*{Acknowledgments}
We thank Graham Denham and Michael Falk for useful conversations.

\newcommand{\arxiv}[1]
{\texttt{\href{http://arxiv.org/abs/#1}{arxiv:#1}}}

\renewcommand{\MR}[1]
{\href{http://www.ams.org/mathscinet-getitem?mr=#1}{MR#1}}

\end{document}